\documentclass[11pt]{article}%
\usepackage{amsfonts}
\usepackage{graphicx}
\usepackage{amsmath}
\usepackage{amssymb}%
\providecommand{\U}[1]{\protect\rule{.1in}{.1in}}
\newtheorem{theorem}{Theorem}

\newtheorem{claim}[theorem]{Claim}

\newtheorem{corollary}[theorem]{Corollary}

\newtheorem{example}[theorem]{Example}

\newtheorem{lemma}[theorem]{Lemma}

\newtheorem{proposition}[theorem]{Proposition}
\newtheorem{remark}[theorem]{Remark}

\begin{document}
Journal of Applied Probability 49, 4 1091-1105\bigskip\ (2012). Applied
Probability Trust

\begin{center}
{\Large On the functional CLT for reversible Markov Chains with nonlinear
growth of the variance}

\bigskip

Martial Longla, Costel Peligrad and Magda Peligrad\footnote{Supported in part
by a Charles Phelps Taft Memorial Fund grant and the NSA\ grant
H98230-11-1-0135, and the NSF grant DMS-1208237.}

\bigskip
\end{center}

Department of Mathematical Sciences, University of Cincinnati, PO Box 210025,
Cincinnati, Oh 45221-0025, USA.

Email: martiala@mail.uc.edu, peligrc@ucmail.uc.edu, peligrm@ucmail.uc.edu

\textit{Key words}: Maximal inequality; reversible processes; Markov chains;
martingale approximation; tightness; functional central limit theorem.

\textit{Mathematical Subject Classification} (2000): 60 F 17, 60 G 05, 60 G
10.\bigskip

\begin{center}
Abstract

\bigskip
\end{center}

In this paper we study the functional central limit theorem for stationary
Markov chains with self-adjoint operator and general state space. We
investigate the case when the variance of the partial sum is not
asymptotically linear in $n,$ and establish that conditional convergence in
distribution of partial sums\ implies functional CLT. The main tools are
maximal inequalities that are further exploited to derive conditions for
tightness and convergence to the Brownian motion.

\section{Introduction}

Kipnis and Varadhan (1986) showed that for an additive functional zero mean
$S_{n}$ of a stationary reversible Markov chain the condition $var(S_{n}%
)/n\rightarrow\sigma^{2}$ implies convergence of $S_{[nt]}/\sqrt{n}$ to the
Brownian motion (here, $[nt]$ is the integer part of $nt$). There is a
considerable amount of papers that further extend and apply this result to
infinite particle systems, random walks, processes in random media,
Metropolis-Hastings algorithms. Among others, Kipnis and Landim (1999)
considered interacting particle systems, Tierney (1994) discussed the
applications to Markov Chain Monte Carlo. Wu (1999) and Zhao and Woodroofe
(2008) studied the law of the iterated logarithm, \ Derriennic and Lin (2001)
and Cuny and Peligrad (2012) investigated the central limit theorem started at
a point.

Recently, Zhao et al. (2010) addressed the conditional central limit theorem
question under the weaker condition $var(S_{n})=nh(n)$, where $h$ is a slowly
varying function (i.e. $\lim_{n\rightarrow\infty}h(nt)/h(n)=1$ for all $t>0$).
They showed by example the surprising result that the distribution of
$S_{[nt]}/stdev(S_{n})$ needs not converge to the standard normal distribution
in this case. They developed sufficient conditions for convergence to a
(possibly non-standard) normal distribution imposed to an approximating martingale.

In this paper we address the question of functional central limit theorem for
the case considered by Zhao et al. (2010). Our goal is to establish sufficient
conditions imposed on the original sequence. We also show that for reversible
Markov chains conditional convergence in distribution of partial sums properly
normalized implies functional CLT. The main tools to prove this result are new
maximal inequalities based on a triangular forward-backward martingale
decomposition and tightness results.

Our paper is organized as follows: Section 2 contains the definitions, a short
background of the problem and the results. Section 3 is devoted to the proofs.
Section 4 contains a functional central limit theorem for an additive
functional associated to a Metropolis-Hastings algorithm, with the variance of
partial sums behaving asymptotically like $nh(n)$ (where $h$ is a slowly
varying function). All throughout the paper $\Rightarrow$ denotes weak
convergence, $[x]$ is the integer part of $x$ and $\rightarrow^{\mathbb{P}}%
$denotes convergence in probability. The notation \thinspace$a_{n}\sim b_{n}$
means $a_{n}/b_{n}\rightarrow1$ as $n\rightarrow\infty;$ $a_{n}=o(b_{n})$
means $a_{n}/b_{n}\rightarrow0$ as $n\rightarrow\infty.$

\section{Definitions, background and results}

We assume that $(\xi_{n})_{n\in\mathbb{Z}}$ is a stationary Markov chain
defined on a probability space $(\Omega,\mathcal{F},\mathbb{P})$ with values
in a general state space $(S,\mathcal{A})$. The marginal distribution is
denoted by $\pi(A)=\mathbb{P}(\xi_{0}\in A)$. Assume there is a regular
conditional distribution for $\xi_{1}$ given $\xi_{0}$ denoted by
$Q(x,A)=\mathbb{P}(\xi_{1}\in A|\,\xi_{0}=x)$. Let $Q$ also denote the Markov
operator {acting via $(Qf)(x)=\int_{S}f(s)Q(x,ds).$ Next, let $\mathbb{L}%
_{0}^{2}(\pi)$ be the set of measurable functions on $S$ such that $\int
f^{2}d\pi<\infty$ and $\int fd\pi=0.$ For some function } ${f}\in$%
{$\mathbb{L}_{0}^{2}(\pi)$, let }%
\begin{equation}
{X_{i}=f(\xi_{i}),\ S_{n}=\sum\limits_{i=1}^{n}X_{i},\ }\sigma_{n}%
=(\mathbb{E}S_{n}^{2})^{1/2}. \label{defX}%
\end{equation}
{\ Denote by $\mathcal{F}_{k}$ the $\sigma$--field generated by $\xi_{i}$ with
$i\leq k$. }

For any integrable random variable $X$ we denote $\mathbb{E}_{k}%
(X)=\mathbb{E}(X|\mathcal{F}_{k}).$ Under this notation, $\mathbb{E}_{0}%
(X_{1})=(Qf)(\xi_{0})=\mathbb{E}(X_{1}|\xi_{0}).$ We denote by ${{||X||}_{p}}$
the norm in {$\mathbb{L}_{p}$}$(\Omega,\mathcal{F},\mathbb{P}).$ {\ }

The Markov chain is called reversible if $Q=Q^{\ast},$ where $Q^{\ast}$ is the
adjoint operator of \ $Q$. The condition of reversibility is equivalent to
requiring that $(\xi_{0},\xi_{1})$ and $(\xi_{1},\xi_{0})$ have the same
distribution or%
\[
\int_{A}Q(\omega,B)\pi(d\omega)=\int_{B}Q(\omega,A)\pi(d\omega)
\]
for all Borel sets $A,B\in\mathcal{A}$.

Kipnis and Varadhan (1986) assumed that%
\begin{equation}
\lim_{n\rightarrow\infty}\frac{\sigma_{n}^{2}}{n}=\sigma_{f}^{2}\text{ }
\label{var}%
\end{equation}
and proved that for any reversible Markov chain defined by (\ref{defX}) this
condition implies
\begin{equation}
\frac{S_{[nt]}}{\sqrt{n}}\Rightarrow|\sigma_{f}|W(t)\text{,} \label{IP}%
\end{equation}
where $W(t)$ is the standard Brownian motion.

Recently Zhao et al. (2010) analyzed the case when%
\begin{equation}
\sigma_{n}^{2}=nh(n),\text{ with }h\ \text{a slowly varying function.}
\label{cond sigma}%
\end{equation}

In their Proposition 1, they showed that without loss of generality, one can
assume that $h(n)\rightarrow\infty,$ since otherwise either (\ref{var}) holds,
(and this case is already known) or $2S_{n}=(1+(-1)^{n-1})${$X_{1\ }$ $a.s$}.
Then, in their Proposition 2 they showed that the representation
(\ref{cond sigma}) implies
\begin{equation}
||\mathbb{E}_{0}(S_{n})||_{2}=o(\sigma_{n}). \label{condMA}%
\end{equation}
On the other hand it is well known that (\ref{condMA}) implies
(\ref{cond sigma}); see for instance Lemma 1 in Wu and Woodroofe (2004).
Therefore, we can state Proposition 2 in Zhao et al. (2010) as follows:

\begin{proposition}
\label{equivalent}For a stationary reversible Markov chain $(X_{n}%
)_{n\in\mathbb{Z}}$ defined by (\ref{defX}), the relations (\ref{cond sigma})
and (\ref{condMA}) are equivalent.
\end{proposition}

In their Corollary 2, Zhao et al. (2010) gave sufficient conditions for the
validity of the conditional CLT in terms of conditions imposed on the
differences of an approximating martingale. In addition, they provided an
example of reversible Markov chain satisfying (\ref{cond sigma}), for which
the central limit theorem holds with a different normalization.

Throughout this paper we shall assume that $\sigma_{n}^{2}\rightarrow\infty.$

By conditional convergence in distribution, denoted by $Y_{n}|\mathcal{F}%
_{0}\Rightarrow Y,$ we understand that for any function $g$ which is
continuous and bounded
\[
\mathbb{E}_{0}(g(Y_{n}))\rightarrow^{\mathbb{P}}\mathbb{E}g(Y)\text{ as
}n\rightarrow\infty.
\]
In other words, let $\mathbb{P}^{x}$ be the probability associated with the
Markov chain started from $x$ and let $\mathbb{E}^{x}$ be the corresponding
expectation. Then, for any $\varepsilon>0$
\[
\pi\{x:|\mathbb{E}^{x}g(Y_{n})-\mathbb{E}g(Y)|>\varepsilon\}\rightarrow
0\text{.}%
\]

One of our results is the following invariance principle for functionals of
stationary reversible Markov chains. Define%
\[
W_{n}(t)=\frac{S_{[nt]}}{\sigma_{n}}%
\]

\begin{theorem}
\label{IP1} Assume $(\xi_{n})_{n\in\mathbb{Z}}$ is a stationary reversible
Markov chain as defined above. Define $(X_{i})_{i\in\mathbb{Z}}$ by
(\ref{defX}) and assume that (\ref{cond sigma}) is satisfied and $S_{n}%
/\sigma_{n}$ is conditionally convergent in distribution to $L$. Then,
\begin{equation}
\text{\ \ \ \ \ \ \ \ \ \ \ \ \ }W_{n}(t)\Rightarrow cW(t),\text{
\ \ \ \ \ \ \ \ \ \ \ \ \ \ \ \ \ \ \ \ \ \ \ } \label{IP2}%
\end{equation}
where $W(t)$ is a standard Brownian motion and $c$ is the standard deviation
of $L$.
\end{theorem}

Theorem \ref{IP1} does not require special properties of the Markov chain such
as irreducibility and aperiodicity. However, if these properties are satisfied
we have the following simplification:

\begin{corollary}
\label{corIA} Assume $(\xi_{n})_{n\in\mathbb{Z}}$ is a stationary, reversible,
irreducible and aperiodic Markov chain such that (\ref{cond sigma}) is
satisfied. Then $S_{n}/\sigma_{n}\Rightarrow L$ implies\ (\ref{IP2}).
\end{corollary}

The proof of Theorem \ref{IP1} requires the development of several tools.
First, we shall establish maximal inequalities that have interest in
themselves. As in the Doob maximal inequalities for martingales, we shall
compare moments and tail distributions of the maximum of partial sums with
those of the corresponding partial sums.

\begin{proposition}
\label{maxrev}Let $(X_{i})_{i\in\mathbb{Z}}$ be defined by (\ref{defX}) and
$Q=Q^{\ast}$. Let $p>1$ and $q>1$ such that $1/p+1/q=1$. Then for all
$n\geq1,$%
\begin{equation}
||\max_{1\leq i\leq n}|S_{i}|\text{ }||_{p}\leq||\max_{1\leq i\leq n}%
|X_{i}|\text{ }||_{p}+(4q+3)\max_{1\leq i\leq n}||S_{i}||_{p}\text{.}
\label{maxL2}%
\end{equation}

\end{proposition}

\begin{remark}
Let $p=2$. Since $(X_{i})_{i\in\mathbb{Z}}$ is stationary it is well known
that
\[
||\max_{1\leq i\leq n}|X_{i}|\text{ }||_{2}=o(n^{1/2})\text{ as }%
n\rightarrow\infty.
\]
If we assume in addition $\lim\inf_{n}\sigma_{n}^{2}/n>0$ we deduce that there
exists $C>0$, such that
\[
||\max_{1\leq i\leq n}|S_{i}|\text{ }||_{2}\leq C\max_{1\leq i\leq n}%
||S_{i}||_{2}\text{.}%
\]

\end{remark}

For the proof of tightness, it is also convenient to have inequalities for the
tail probabilities of partial sums. We shall also establish:

\begin{proposition}
\label{maxprob} Let $(X_{i})_{i\in\mathbb{Z}}$ be defined by (\ref{defX}) and
$Q=Q^{\ast}$. Then, for every $x>0$ and $n\geq1,$
\[
\mathbb{P}(\max_{1\leq i\leq n}|S_{i}|>x)\leq\frac{2}{x}[18\mathbb{E}%
|S_{n}|I(|S_{n}|>x/12)+55\max_{1\leq i\leq n}||\mathbb{E}_{0}(S_{i}%
)||_{1}+||\max_{1\leq i\leq n}|X_{i}|\text{ }||_{1}].
\]
\textbf{\ }
\end{proposition}

An important step in the proof of Theorem \ref{IP1} is the use of tightness
conditions. We shall give two necessary conditions for tightness, that will
ensure continuity of every limiting process.

\begin{proposition}
\label{tight} Assume $X_{i}$ is defined by (\ref{defX}), condition
(\ref{cond sigma}) is satisfied and one of the following two conditions holds:

\begin{enumerate}
\item $(S_{n}^{2}/\sigma_{n}^{2})_{n\geq1}$ is uniformly integrable;

\item $S_{n}/\sigma_{n}$ is convergent in distribution.
\end{enumerate}

Then, $W_{n}(t)$ is tight in $D(0,1)$ endowed with uniform topology and any
limiting process is continuous.
\end{proposition}

Finally, we give sufficient conditions for convergence to the standard
Brownian Motion.

\begin{proposition}%
\index{IP}%
\label{necsuf} Assume $(\xi_{n})_{n\in\mathbb{Z}}$ is a stationary reversible
Markov chain. Define $(X_{i})_{i\in\mathbb{Z}}$ by (\ref{defX}) and assume
that (\ref{condMA}) is satisfied. Assume $(S_{n}^{2}/\sigma_{n}^{2})_{n\geq1}$
is uniformly integrable and
\begin{equation}
\lim_{n\rightarrow\infty}\frac{||\mathbb{E}_{0}(S_{n}^{2})-\sigma_{n}%
^{2}||_{1}}{\sigma_{n}^{2}}=0\text{.} \label{key2}%
\end{equation}
Then,
\[
W_{n}(t)\Rightarrow W(t)\text{.}%
\]

\end{proposition}

\section{Proofs}

We start with a preliminary martingale decomposition that combines ideas from
Wu and Woodroofe (2004) with forward-backward martingale approximation of
Meyer and Zheng (1984) and Lyons and Zheng (1988).

\subsection{Forward-backward martingale decomposition}

As in Wu and Woodroofe (2004) for $n\geq1$ fixed, define the stationary
sequences
\[
\theta_{k}^{n}=\frac{1}{n}\sum_{i=0}^{n-1}\mathbb{E}_{k}(X_{k}+...+X_{k+i}%
)\text{, and}%
\]%
\begin{equation}
D_{k}^{n}=\theta_{k}^{n}-\mathbb{E}_{k-1}(\theta_{k}^{n})\text{.}
\label{martdef}%
\end{equation}
Then, $(D_{k}^{n})_{k\in Z}$ is a triangular array of martingale differences
adapted to the filtration $\mathcal{F}_{n}=\sigma(\xi_{i},i\leq n)$. Notice
that%
\begin{align*}
\theta_{k}^{n}  &  =X_{k}+\frac{1}{n}\sum_{i=1}^{n-1}\mathbb{E}_{k}%
(S_{k+i}-S_{k})=X_{k}+\mathbb{E}_{k}(\theta_{k+1}^{n})-\frac{1}{n}%
\mathbb{E}_{k}(S_{k+n}-S_{k})\\
&  =X_{k}+\theta_{k+1}^{n}-D_{k+1}^{n}-\frac{1}{n}\mathbb{E}_{k}(S_{k+n}%
-S_{k}).
\end{align*}
Therefore,
\begin{equation}
X_{k}=D_{k+1}^{n}+\theta_{k}^{n}-\theta_{k+1}^{n}+\frac{1}{n}\mathbb{E}%
_{k}(S_{k+n}-S_{k}). \label{decX}%
\end{equation}
We construct now a martingale approximation for the reversed process adapted
to the filtration $\mathcal{G}_{n}=\sigma(\xi_{i},i\geq n)$. We introduce the
notation $\mathbb{\tilde{E}}_{1}(X_{0})=\mathbb{E}(X_{0}|\mathcal{G}%
_{1})\ =\mathbb{E}(X_{0}|\xi_{1})=(Q^{\ast}f)(\xi_{1}).$

Now, let%
\[
\tilde{\theta}_{k}^{n}=\frac{1}{n}\sum_{i=0}^{n-1}\mathbb{\tilde{E}}%
_{k}(X_{k-i}+...+X_{k}).
\]
With this notation
\begin{equation}
X_{k+1}=\tilde{D}_{k}^{n}+\tilde{\theta}_{k+1}^{n}-\tilde{\theta}_{k}%
^{n}+\frac{1}{n}\mathbb{\tilde{E}}_{k+1}(X_{-n+k+1}+...+X_{k})\text{.}
\label{defXr}%
\end{equation}
where $\tilde{D}_{k}^{n}$ are martingale differences with respect to the
filtration $\mathcal{G}_{k}=\sigma(\xi_{i},i\geq k),$ $\tilde{D}_{k}%
^{n}=\tilde{\theta}_{k}^{n}-\mathbb{E}_{k+1}\tilde{\theta}_{k}^{n}$.

If we assume that $Q=Q^{\ast}$, we have $\mathbb{\tilde{E}}_{1}(X_{0})=$
$\mathbb{E}(X_{2}|\xi_{1})=(Qf)(\xi_{1})$. Therefore, $\tilde{\theta}_{k}%
^{n}=\theta_{k}^{n}$, $\tilde{\theta}_{k+1}^{n}=\theta_{k+1}^{n}$ and
$\mathbb{\tilde{E}}_{k+1}(X_{-n+k+1}+...+X_{k})=\mathbb{E}_{k+1}%
(X_{k+2}+...+X_{k+n+1})$. Adding relations (\ref{decX}) and (\ref{defXr})
leads to
\[
X_{k}+X_{k+1}=D_{k+1}^{n}+\tilde{D}_{k}^{n}+\frac{1}{n}\mathbb{E}_{k}%
(S_{n}-S_{k})+\frac{1}{n}\mathbb{E}_{k+1}(S_{k+n+1}-S_{k+1}).
\]
Summing these relations we obtain the representation
\[
\sum_{i=0}^{k-1}(X_{i}+X_{i+1})=\sum_{i=1}^{k}[(D_{i}^{n}+\tilde{D}_{i-1}%
^{n})+\frac{1}{n}\mathbb{E}_{i-1}(S_{n+i-1}-S_{i-1})+\frac{1}{n}\mathbb{E}%
_{i}(S_{n+i}-S_{i})].
\]
So,%
\[
2S_{k}+(X_{0}-X_{k})=\sum_{i=1}^{k}(D_{i}^{n}+\tilde{D}_{i-1}^{n})+\bar{R}%
_{k}^{n},
\]
where%
\[
\bar{R}_{k}^{n}=\frac{1}{n}\sum_{i=1}^{k}[\mathbb{E}_{i-1}(S_{n+i-1}%
-S_{i-1})+\mathbb{E}_{i}(S_{n+i}-S_{i})]\text{.}%
\]
Therefore, in the reversible case, we get the following forward-backward
martingale representation%
\begin{equation}
S_{k}=\frac{1}{2}[(X_{k}-X_{0})+(M_{k}^{n}+\tilde{M}_{k}^{n})+\bar{R}_{k}%
^{n}],\label{FB-martapp}%
\end{equation}
where $M_{k}^{n}=\sum_{i=1}^{k}D_{i}^{n}$ is a forward martingale adapted to
the filtration {$\mathcal{F}_{k}$ and }$\tilde{M}_{k}^{n}=\sum_{i=0}%
^{k-1}\tilde{D}_{i}^{n}$ is a backward martingale adapted to the filtration
{$\mathcal{G}_{k}.$}

Also, it is convenient to point out a related martingale approximation, which
helps us relate the partial sums with a martingale adapted to the same
filtration. Notice that
\[
\theta_{k}^{n}=X_{k}+\frac{1}{n}\sum_{i=1}^{n-1}\mathbb{E}_{k}(S_{k+i}%
-S_{k})=X_{k}+\bar{\theta}_{k}^{n},
\]%
\[
\text{where}\quad\ \ \bar{\theta}_{k}^{n}=\frac{1}{n}\sum_{i=1}^{n-1}%
\mathbb{E}_{k}(S_{k+i}-S_{k}).
\]
Starting from (\ref{decX}) and using this notation we obtain
\[
X_{k+1}=D_{k+1}^{n}+\bar{\theta}_{k}^{n}-\bar{\theta}_{k+1}^{n}+\frac{1}%
{n}\mathbb{E}_{k}(S_{k+n}-S_{k}).
\]
So, summing these relations, denoting as above $M_{k}^{n}=\sum_{i=1}^{k}%
D_{i}^{n}$, we obtain for every stationary sequence, not necessarily
reversible, and for any $n$ and $m,$%
\begin{gather}
S_{m}=M_{m}^{n}+R_{m}^{n},\label{martdec}\\
\text{where \ \ }R_{m}^{n}=\bar{\theta}_{0}^{n}-\bar{\theta}_{m}^{n}+\frac
{1}{n}\sum_{k=0}^{m-1}\mathbb{E}_{k}(S_{k+n}-S_{k}).\nonumber
\end{gather}

\subsection{Proof of Proposition \ref{maxrev}}

We start from (\ref{FB-martapp}) and take the maximum on both sides. We easily obtain%

\begin{equation}
\max_{1\leq i\leq n}|S_{i}|\leq\frac{1}{2}(|X_{0}|+\max_{1\leq i\leq n}%
|X_{i}|+\max_{1\leq i\leq n}|M_{i}^{n}+\tilde{M}_{i}^{n}|+\max_{1\leq i\leq
n}|\bar{R}_{i}^{n}|)\text{.} \label{martdecmax}%
\end{equation}
Notice that%
\[
\ \max_{1\leq i\leq n}|\bar{R}_{i}^{n}|\leq\frac{1}{n}\sum_{i=1}%
^{n}(|\mathbb{E}_{i-1}(S_{n+i-1}-S_{i-1})|+|\mathbb{E}_{i}(S_{n+i}%
-S_{i})|)\text{,}%
\]
whence, by Minkowski's inequality and stationarity, for any $p\geq1$%
\begin{equation}
||\max_{1\leq i\leq n}|\bar{R}_{i}^{n}|\text{ }||_{p}\leq\frac{2}{n}\sum
_{i=1}^{n}||\mathbb{E}_{i}(S_{n+i}-S_{i})||_{p}=2||\mathbb{E}_{0}(S_{n}%
)||_{p}\text{.} \label{errorL_p}%
\end{equation}
Taking into account that $\max_{1\leq k\leq n}|\tilde{M}_{k}^{n}|\leq
|\tilde{M}_{n}^{n}|+\max_{1\leq k\leq n}|\tilde{M}_{n}^{n}-\tilde{M}_{k}%
^{n}|,$ we easily deduce
\[
||\max_{1\leq k\leq n}|M_{k}^{n}+\tilde{M}_{k}^{n}|\text{ }||_{p}\leq
||\max_{1\leq k\leq n}|M_{k}^{n}|\text{ }||_{p}+||\max_{1\leq k\leq n}%
|\tilde{M}_{n}^{n}-\tilde{M}_{k}^{n}|\text{ }||_{p}+||\tilde{M}_{n}^{n}%
||_{p},
\]
whence, by Doob maximal inequality applied twice, stationarity and
reversibility,%
\[
||\max_{1\leq k\leq n}|M_{k}^{n}+\tilde{M}_{k}^{n}|\text{ }||_{p}\leq
q||M_{n}^{n}||_{p}+(q+1)||\tilde{M}_{n}^{n}||_{p}=(2q+1)||M_{n}^{n}||_{p},
\]
(where $q$ is the conjugate of $p)$.

From (\ref{martdec}) we have $M_{n}^{n}=S_{n}-R_{n}^{n}$, and from Minkowski's
inequality we deduce that%
\[
||M_{n}^{n}||_{p}\leq||S_{n}||_{p}+\frac{2}{n}\sum\limits_{i=0}^{n-1}%
||\mathbb{E}_{0}(S_{i})||_{p}+||\mathbb{E}_{0}(S_{n})||_{p}%
\]
whence,%
\begin{equation}
||M_{n}^{n}||_{p}\leq||S_{n}||_{p}+3\max_{1\leq i\leq n}||\mathbb{E}_{0}%
(S_{i})||_{p}\text{.} \label{estmart}%
\end{equation}
From (\ref{martdecmax}), (\ref{errorL_p}) and (\ref{estmart})\ we deduce the
following extension of the Doob maximal inequality for reversible processes:
\[
||\max_{1\leq i\leq n}|S_{i}|\text{ }||_{p}\leq\frac{1}{2}(||X_{0}%
||_{p}+||\max_{1\leq i\leq n}|X_{i}|\text{ }||_{p}+(2q+1)[||S_{n}||_{p}%
+3\max_{1\leq i\leq n}||\mathbb{E}_{0}(S_{i})||_{p}]+2||\mathbb{E}_{0}%
(S_{n})||_{p}).
\]
Taking now into account that $||\mathbb{E}_{0}(S_{i})||_{p}\leq||S_{i}||_{p},$
Proposition \ref{maxrev} is established. $\Diamond$

\subsection{\textbf{\ Proof of Proposition \ref{maxprob}}}

For the proof of this proposition we shall use the following claim that can be
easily obtained by truncation:

\begin{claim}
Let $X$ and $Y$ be two positive random variables. Then for all $x\geq0$
\[
\mathbb{E}XI(Y>x)\leq\mathbb{E}XI(X>x/2)+\frac{x}{2}\mathbb{P}(Y>x).
\]

\end{claim}

For every $x\geq0$, using (\ref{martdecmax}), we obtain%
\begin{equation}
\mathbb{P}(\max_{1\leq i\leq n}|S_{i}|>x)\leq\mathbb{P}(\max_{1\leq i\leq
n}|M_{i}^{n}+\tilde{M}_{i}^{n}|>x)+\mathbb{P}(|X_{0}|+\max_{1\leq i\leq
n}|X_{i}|+\max_{1\leq i\leq n}|\bar{R}_{i}^{n}|>x). \label{3a}%
\end{equation}
Applying the Markov inequality, then the triangle inequality followed by
(\ref{errorL_p}) with $p=1,$ leads to
\begin{equation}
\mathbb{P}(|X_{0}|+\max_{1\leq i\leq n}|X_{i}|+\max_{1\leq i\leq n}|\bar
{R}_{i}^{n}|>x)\leq\frac{2}{x}(||\max_{1\leq i\leq n}|X_{i}|\text{ }%
||_{1}+||\mathbb{E}_{0}(S_{n})||_{1}). \label{3b}%
\end{equation}
By triangle inequality and reversibility%
\begin{gather*}
\mathbb{P}(\max_{1\leq i\leq n}|M_{i}^{n}+\tilde{M}_{i}^{n}|>x)\leq
\mathbb{P}(\max_{1\leq i\leq n}|M_{i}^{n}|>x/3)+\\
\ +\mathbb{P}(\max_{1\leq i\leq n}|\sum_{k=i}^{n}\tilde{D}_{i}^{n}%
|>x/3)+\mathbb{P}(|\tilde{M}_{n}^{n}|>x/3)\leq3\mathbb{P}(\max_{1\leq i\leq
n}|M_{i}^{n}|>x/3).
\end{gather*}
Then, by Doob maximal inequality and the above Claim applied to $X=|M_{n}%
^{n}|$ and $Y=\max_{1\leq i\leq n}|M_{i}^{n}|$ we obtain%
\begin{gather*}
\mathbb{P}(\max_{1\leq i\leq n}|M_{i}^{n}|>x/3)\leq\frac{3}{x}\mathbb{E}%
|M_{n}^{n}|I(\max_{1\leq i\leq n}|M_{i}^{n}|>x/3)\\
\leq\frac{3}{x}\mathbb{E}|M_{n}^{n}|I(|M_{n}^{n}|>x/6)+\frac{1}{2}%
\mathbb{P}(\max_{1\leq i\leq n}|M_{i}^{n}|>x/3),
\end{gather*}
implying
\[
\mathbb{P}(\max_{1\leq i\leq n}|M_{i}^{n}|>x/3)\leq\frac{6}{x}\mathbb{E}%
|M_{n}^{n}|I(|M_{n}^{n}|>x/6)\text{.}%
\]
Now, we express the right-hand side in terms of $S_{n}.$ By (\ref{martdec}) we
have $M_{n}^{n}=S_{n}-R_{n}^{n}$ and using the fact that for all positive real
numbers $x,y,a$ we have $(x+y)I(x+y>a)\leq2xI(x>a/2)+2yI(y>a/2)\leq$
$2xI(x>a/2)+2y,$ we obtain%
\begin{gather*}
\mathbb{E}|M_{n}^{n}|I(|M_{n}^{n}|>x/6)\leq2\mathbb{E}|S_{n}|I(|S_{n}%
|>x/12)+2||R_{n}^{n}||_{1}\\
\leq2\mathbb{E}|S_{n}|I(|S_{n}|>x/12)+6\max_{1\leq i\leq n}||\mathbb{E}%
_{0}(S_{i})||_{1}\text{.}%
\end{gather*}
Therefore,%
\[
\mathbb{P}(\max_{1\leq i\leq n}|M_{i}^{n}|>x/3)\leq\frac{6}{x}[2\mathbb{E}%
|S_{n}|I(|S_{n}|>x/12)+6\max_{1\leq i\leq n}||\mathbb{E}_{0}(S_{i})||_{1}]
\]
and so%
\begin{equation}
\mathbb{P}(\max_{1\leq k\leq n}|M_{k}^{n}+\tilde{M}_{k}^{n}|>x)\leq\frac
{18}{x}[2\mathbb{E}|S_{n}|I(|S_{n}|>x/12)+6\max_{1\leq i\leq n}||\mathbb{E}%
_{0}(S_{i})||_{1}]\text{.} \label{3c}%
\end{equation}
Thus, (\ref{3a}), (\ref{3b}) and (\ref{3c}) lead to
\[
\mathbb{P}(\max_{1\leq i\leq n}|S_{i}|>x)\leq\frac{2}{x}[18\mathbb{E}%
|S_{n}|I(|S_{n}|>x/12)+55\max_{1\leq i\leq n}||\mathbb{E}_{0}(S_{i}%
)||_{1}+||\max_{1\leq i\leq n}|X_{i}|\text{ }||_{1}].
\]
$\Diamond$

\subsection{\textbf{Proof of Proposition \ref{tight}}}

We prove first the conclusion of the proposition under the assumption that
$(S_{n}^{2}/\sigma_{n}^{2})_{n\geq1}$ is uniformly integrable.

By stationarity and by Theorem 8.3 in Billingsley (1968) formulated for random
elements of D (see page 137 in Billingsley, 1968) we have to show that for all
$\varepsilon>0$
\begin{equation}
\lim_{\delta\rightarrow0^{+}}\lim\sup_{n\rightarrow\infty}\frac{1}{\delta
}\mathbb{P}(\max_{1\leq k\leq\lbrack n\delta]}|S_{k}|>\varepsilon\sigma
_{n})=0\text{.} \label{tight to show}%
\end{equation}
By Proposition \ref{maxprob},
\begin{align}
\mathbb{P}(\max_{1\leq k\leq\lbrack n\delta]}|S_{k}|  &  >\varepsilon
\sigma_{n})\leq\frac{2}{\varepsilon\sigma_{n}}[18\mathbb{E}|S_{[n\delta
]}|I(|S_{[n\delta]}|>\varepsilon\sigma_{n}/12)+\label{rel1}\\
&  55\max_{1\leq i\leq\lbrack n\delta]}\mathbb{E}|\mathbb{E}_{0}%
(S_{i})|+\mathbb{E}\max_{1\leq i\leq n}|X_{i}|]\text{.}\nonumber
\end{align}
We shall analyze each term from the right-hand side of inequality (\ref{rel1}) separately.

By the fact that $\ \lim_{n\rightarrow\infty}\sigma_{\lbrack n\delta]}%
^{2}/\delta\sigma_{n}^{2}=1$, taking into account uniform integrability of
$(S_{n}^{2}/\sigma_{n}^{2})_{n\geq1}$ leads to
\begin{align*}
\lim_{\delta\rightarrow0^{+}}\lim\sup_{n\rightarrow\infty}\frac{1}%
{\delta\sigma_{n}}\mathbb{E}|S_{[n\delta]}|I(|S_{[n\delta]}|  &
>\frac{\varepsilon\sigma_{n}}{12})\leq\\
\lim_{\delta\rightarrow0^{+}}\lim\sup_{n\rightarrow\infty}\frac{24}%
{\varepsilon\sigma_{n}^{2}}\mathbb{E}S_{n}^{2}I(\frac{|S_{n}|}{\sigma_{n}}  &
>\frac{\varepsilon}{24\delta^{1/2}})=0.
\end{align*}
By stationarity and the fact that $\lim\inf_{n}\sigma_{n}^{2}/n>0$ we have
\begin{equation}
\frac{1\text{ }}{\sigma_{n}^{2}}(\mathbb{E}\max_{1\leq i\leq n}|X_{i}%
|)^{2}\leq\frac{1\text{ }}{\sigma_{n}^{2}}\mathbb{E}\max_{1\leq i\leq n}%
|X_{i}|^{2}\rightarrow0\text{ as }n\rightarrow\infty. \label{conv max X}%
\end{equation}
Then, by condition (\ref{condMA}) and Proposition \ref{equivalent},
\begin{equation}
\frac{1}{\sigma_{n}^{2}}\max_{1\leq i\leq\lbrack n\delta]}(\mathbb{E}%
|\mathbb{E}_{0}(S_{i})|)^{2}\leq\frac{1}{\sigma_{n}^{2}}\max_{1\leq
i\leq\lbrack n\delta]}\mathbb{E}[\mathbb{E}_{0}(S_{i})]^{2}\rightarrow0\text{
as }n\rightarrow\infty. \label{second term}%
\end{equation}
Then, combining the last three convergence results with the inequality
(\ref{rel1}) leads to (\ref{tight to show}).

To prove the second part of this proposition, assume now that $S_{n}%
/\sigma_{n}\Rightarrow L.$ By Theorem 5.3 in Billingsley (1968), we notice
that the limit has finite second moment, namely%
\begin{equation}
\mathbb{E}L^{2}\leq\lim\inf_{n\rightarrow\infty}||S_{n}||_{2}^{2}/\sigma
_{n}^{2}=1. \label{moment 2}%
\end{equation}
Furthermore, since $(|S_{n}|/\sigma_{n})_{n\geq1}$ is uniformly integrable
(because $\mathbb{E}S_{n}^{2}/\sigma_{n}^{2}=1$), by (\ref{condMA}) and
Theorem 5.4 in Billingsley (1968), it follows that%
\begin{align}
\lim\sup_{n\rightarrow\infty}\frac{1}{\sigma_{n}}\mathbb{E}|S_{[n\delta
]}|I(|S_{[n\delta]}|  &  >\frac{\varepsilon\sigma_{n}}{12})\leq\frac{1}%
{\sqrt{\delta}}\lim_{n\rightarrow\infty}\frac{1}{\sigma_{\lbrack n\delta]}%
}\mathbb{E}|S_{[n\delta]}|I(\frac{|S_{[n\delta]}|}{\sigma_{\lbrack n\delta]}%
}>\frac{\varepsilon}{24\sqrt{\delta}})\label{limB}\\
&  =\frac{1}{\sqrt{\delta}}\mathbb{E}|L|I(|L|>\frac{\varepsilon}%
{24\sqrt{\delta}}).\nonumber
\end{align}
By passing to the limit in relation (\ref{rel1}) and using (\ref{conv max X}),
(\ref{second term}), and (\ref{limB}) we obtain,
\[
\lim\sup_{n\rightarrow\infty}\frac{1}{\delta}\mathbb{P}(\max_{1\leq
k\leq\lbrack n\delta]}|S_{k}|>\varepsilon\sigma_{n})\leq\frac{36}%
{\varepsilon\delta^{1/2}}\mathbb{E}|L|I(|L|>\frac{\varepsilon}{24\sqrt{\delta
}})\text{.}%
\]
Then, clearly
\[
\lim\sup_{n\rightarrow\infty}\frac{1}{\delta}\mathbb{P}(\max_{1\leq
k\leq\lbrack n\delta]}|S_{k}|>\varepsilon\sigma_{n})\leq\frac{36\times
24}{\varepsilon}\mathbb{E}L^{2}I(|L|>\frac{\varepsilon}{24\sqrt{\delta}}).
\]
Finally, taking into account (\ref{moment 2}), the conclusion follows by
letting $\delta\rightarrow0^{+}.$ $\Diamond$

\subsection{\textbf{Proof of Theorem \ref{IP1}}}

Because conditional convergence in distribution implies weak convergence, it
follows that $S_{n}/\sigma_{n}\Rightarrow L$. Then, by the second part of
Proposition \ref{tight}, $W_{n}(t)$ is tight in $C(0,1)$ endowed with uniform
topology with all possible limits in $C(0,1)$. Now, let us consider a
convergent subsequence, say $W_{n^{\prime}}(t)\Rightarrow X(t).$ Then $X(t)$
is continuous and since $S_{n}/\sigma_{n}$ is conditionally convergent in
distribution, $X(t)$ has independent increments (by the next lemma in this
subsection applied on subsequences). It is well known [see, for instance, Doob
(1953), Ch. VIII)] that the process $X(t)$ has the representation
$X(t)=at+bW(t)$ for some constants $a$ and $b,$ where $W(t)$ is the standard
Brownian motion. Without restricting the generality, by symmetry, we can
assume $b>0$. To identify the constants, we use the convergence of moments in
the limit theorem, namely Theorem 5.4 in Billingsley (1968). Notice that
$(S_{n}/\sigma_{n})_{n\geq1}$ is uniformly integrable in $\mathbb{L}_{1}$
since it is bounded in $\mathbb{L}_{2}.$ We use this remark to obtain
\[
\mathbb{E}L=\lim_{n\rightarrow\infty}\mathbb{E}S_{n}/\sigma_{n}=0=\lim
_{n^{\prime}\rightarrow\infty}\mathbb{E}W_{n^{\prime}}(1)=\mathbb{E}%
X(1)=a+b\mathbb{E}W(1)=a,
\]
so $a=0$. Finally, by the same argument it follows that
\[
\lim_{n\rightarrow\infty}\mathbb{E}|S_{n}|/\sigma_{n}=\mathbb{E}%
|L|=\lim_{n^{\prime}\rightarrow\infty}\mathbb{E}|W_{n^{\prime}}%
(1)|=\ \mathbb{E}|X(1)|=b\mathbb{E}|W(1)|=b\sqrt{2/\pi}.
\]
and so $b=\mathbb{E}|L|\sqrt{\pi/2}.$ It follows that $X(t)=(\mathbb{E}%
|L|\sqrt{\pi/2})W(t).$ In particular it follows that $L$ $\ $has normal
distribution and therefore $\mathbb{E}|L|\sqrt{\pi/2}$ is the standard
deviation of $L$. $\Diamond$

\begin{lemma}
\label{LemmaII}Under the assumptions of Theorem \ref{IP1}, if $W_{n}%
(t)\Rightarrow X(t),$ then $X(t)$ has independent increments.
\end{lemma}

Proof. Without loss of generality, for simplicity we consider only two
increments. For any $0\leq s<t\leq1,$ we shall show that%
\[
(W_{n}(s),W_{n}(t)-W_{n}(s))\Rightarrow(X(s),X(t)-X(s))
\]
where $X(s)$ and $X(t)-X(s)$ are independent. By the Cram\'{e}r-Wold device it
is enough to show that for any two real numbers $a$ and $b$,
\begin{gather*}
A=\mathbb{E}\exp[iaW_{n}(s)+ib(W_{n}(t)-W_{n}(s))]\\
-\mathbb{E}\exp[iaX(s)]\mathbb{E}\exp[ib(X(t)-X(s))]\rightarrow0.
\end{gather*}
To see this, notice that
\begin{align*}
&  \mathbb{E}\exp[iaW_{n}(s)+ib(W_{n}(t)-W_{n}(s))]\\
&  =\mathbb{E}\exp[iaW_{n}(s)]\mathbb{E}_{[ns]}\exp[ib(W_{n}(t)-W_{n}(s))].
\end{align*}
By adding and substracting $\mathbb{E}\exp[iaW_{n}(s)]\mathbb{E}%
\exp[ib(X(t)-X(s))]$ to $A$, we easily obtain
\begin{gather*}
|A|\leq\mathbb{E}|\mathbb{E}_{[ns]}\exp[ib(W_{n}(t)-W_{n}(s)]-\mathbb{E}%
\exp[ib(X(t)-X(s))]|+\\
+|\mathbb{E}\exp[iaW_{n}(s)]-\mathbb{E}\exp[iaX(s)]|=I+II.
\end{gather*}
Since we assume that $W_{n}(s)\Rightarrow X(s),$ it follows that
$II\rightarrow0.$ Furthermore, by (\ref{cond sigma}), $X(s)$ and $s^{1/2}L$
are identically distributed.

To treat the term $I$, notice that by stationarity and the definition of
$W_{n}(t)$ we have that
\begin{equation}
I=\mathbb{E}|\mathbb{E}_{0}\exp[{{{ib(S_{[nt]-[ns]}/\sigma_{n})}}}%
]-\mathbb{E}\exp[ib(X(t)-X(s))]|. \label{I1}%
\end{equation}
Because we assume that $\sigma_{n}\rightarrow\infty$ we have
\begin{equation}
\frac{1}{\sigma_{n}}\mathbb{E}|S_{[nt]-[ns]}-S_{[n(t-s)]}|\rightarrow0,
\label{I3}%
\end{equation}
which easily implies that for all $b$,
\begin{equation}
\mathbb{E}|\mathbb{E}_{0}\exp[{{{ib(S_{[nt]-[ns]}/\sigma_{n})}}}%
]-\mathbb{E}_{0}\exp[{{{ib(S_{[n(t-s)]}/\sigma_{n})}}}]|\rightarrow0.
\label{I2}%
\end{equation}
Now, since $S_{[n(t-s)]}/\sigma_{n}\Rightarrow X(t-s)$ \ and $S_{[nt]}%
-S_{[ns]}/\sigma_{n}\rightarrow X(t)-X(s),$ we deduce from (\ref{I3}) and
stationarity that $X(t-s)$ and $X(t)-X(s)$ have the same distribution.
Furthermore, by (\ref{cond sigma}), we deduce that $S_{[n(t-s)]}/\sigma_{n}$
is also conditionally convergent in distribution; so, we also have $X(t-s)$ is
distributed as $(t-s)^{1/2}L.$ Whence, by taking also into account (\ref{I1})
and (\ref{I2}) it follows that
\[
\lim\sup_{n\rightarrow\infty}I=\lim\sup_{n\rightarrow\infty}\mathbb{E}%
|\mathbb{E}_{0}\exp[{{{ib(S_{[n(t-s)]}/\sigma_{n})}}}]-\mathbb{E}%
\exp[ib(X(t-s))]|=0,
\]
leading to the conclusion. $\Diamond$

\subsection{Proof of Corollary \ref{corIA}}

The proof of this corollary follows the lines of Theorem \ref{IP1} with the
exception that we replace Lemma \ref{LemmaII} by the following Lemma:

\begin{lemma}
Under the assumptions of Corollary \ref{corIA}, if $W_{n}(t)\Rightarrow X(t),
$ then $X(t)$ has independent increments.
\end{lemma}

Proof. We mention that by the fact that the Markov chain is stationary,
irreducible and aperiodic it follows that it is absolutely regular [see
Theorems 21.5 and its Corollary 21.7 in volume 2 of Bradley (2007)]. It is
well known that an absolutely regular sequence is strong mixing [see the chart
on page 186, volume 1, Bradley (2007)]. This means that $\alpha_{n}\searrow0$
where%
\[
\alpha_{n}=\sup\mathbb{P}(A\cap B)-\mathbb{P}(A)\mathbb{P}(B);
\]
here the supremum is taken over all $A\in\sigma(\xi_{i},$ $i\leq0)$ and
$B\in\sigma(\xi_{i},$ $i\geq n).$ Because we know from the proof of Theorem
\ref{IP1} that the process $X(t)$ is continuous, it is enough to show that for
all $k$ and $0<s_{1}<t_{1}<s_{2}<t_{2}\,<...<s_{k}<t_{k}<1$ the increments
$(X(t_{i})-X(s_{i}))_{1\leq i\leq k}$ are independent. Now, using the
definitions of $\alpha_{n}$ and $W_{n}(t),$ we get by recurrence
\begin{align*}
|\mathbb{P}(\cap_{i=1}^{k}(W_{n}(t_{i}-s_{i}) &  \in A_{i}))-\Pi_{i=1}%
^{k}\mathbb{P}(W_{n}(t_{i}-s_{i})\in A_{i})|\leq\\
\min_{1\leq i\leq k-1}\alpha_{\lbrack n(s_{i+1}-t_{i})]} &  \rightarrow0\text{
as }n\rightarrow\infty,
\end{align*}
for any Borelians $A_{1},...,A_{k}.$ The conclusion follows by passing to the
limit with $n$. $\Diamond$

\subsection{Proof of Proposition \ref{necsuf}}

\textbf{\ }By Proposition (\ref{equivalent}) we know that $\sigma_{n}%
^{2}=nh(n)$ with $h$ a function slowly varying at infinity. Then, by the first
part of Proposition \ref{tight}, $W_{n}(t)$ is tight in $D(0,1).$ It remains
to apply Theorem 19.4 in Billingsley (1968).

\section{Application to a Metropolis-Hastings algorithm}

In this section we analyze a standardized example of a stationary irreducible
and aperiodic Metropolis-Hastings algorithm with uniform marginal
distribution. This type of Markov chain is interesting since it can easily be
transformed into Markov chains with different marginal distributions. We point
out a central limit theorem under a normalization other than the variance of
partial sums. Markov chains of this type are often studied in the literature
from different points of view, as in Doukhan et al (1994), Rio (2000 and
2009), Merlev\`{e}de and Peligrad (2010). The idea of considering
Metropolis-Hastings algorithm in this context comes from Zhao et al. (2010).

Let $E=[-1,1]$. We define now the transition probabilities of a Markov chain
by
\[
Q(x,A)=(1-|x|)\delta_{x}(A)+|x|\upsilon(A)\,,
\]
where $\delta_{x}$ denotes the Dirac measure and $\upsilon$ on $[-1,1]$
satisfies $\ $%
\begin{equation}
\upsilon(dx)=|x|dx. \label{def-niu}%
\end{equation}
Then, there is$\ $a unique invariant measure, the uniform distribution on
$[-1,1],$
\[
\pi(dx)=dx/2,
\]
and the stationary Markov chain $(\xi_{i})_{i}$ with values in $E$ and
transition probability $Q(x,A)$ is reversible and positively recurrent.
Moreover, for any odd function $f$ \ we have
\begin{equation}
Q^{k}(f)(\xi_{0})=\mathbb{E}(f(\xi_{k})|\xi_{0})=(1-|\xi_{0}|)^{k}f(\xi
_{0})\text{ a.s.} \label{operator}%
\end{equation}
For the odd function $f(x)=sign$ $x$, define $X_{i}=sign$ $\xi_{i}.$ In this
context we shall show:

\begin{example}
Let $(X_{j})_{j\geq1}$ defined above. Then $\sigma_{n}^{2}/(2n\log
n)\rightarrow1$ and
\[
\frac{1}{\sigma_{n}}\sum_{j=1}^{[nt]}X_{j}\Rightarrow\frac{1}{2^{1/2}}W(t).
\]
where $W(t)$ is the standard Brownian motion.
\end{example}

Proof. For any $m\geq0$ \ we have%
\[
\mathbb{E}(X_{0}X_{m})=\mathbb{E}(f(\xi_{0})Q^{m}(f)(\xi_{0}))=\int
_{E}(1-|x|)^{m}\pi(dx)=1/(m+1).\,
\]
Therefore, by simple computations, we obtain%
\[
\sigma_{n}^{2}\sim2n\log n\text{ as }n\rightarrow\infty\text{.}%
\]
\textbf{\ }Now, to find the limiting distribution of $S_{n}$ properly
normalized, we study the regeneration process. Let
\[
T_{0}=\inf\{i>0:\xi_{i}\neq\xi_{0}\}
\]
and
\[
T_{k+1}=\inf\{i>T_{k}:\xi_{i}\neq\xi_{i-1}\},\text{ \ }\tau_{k}=T_{k+1}%
-T_{k}\text{.}%
\]
It is well known that $(\xi_{\tau_{k}},\tau_{k})_{k\geq1}$ are i.i.d. random
variables with $\xi_{\tau_{k}}$ having the distribution $\upsilon.$
Furthermore,
\[
\mathbb{P}(\tau_{1}>n|\xi_{\tau_{1}}=x)=(1-|x|)^{n}\text{.}%
\]
Then, it follows that%
\[
\mathbb{E}(\tau_{1}|\xi_{\tau_{1}}=x)=\frac{1}{|x|}\text{ }\ \ \text{and
\ \ }\mathbb{E}(\tau_{1})=2.
\]
So, by the law of large numbers $T_{n}/n\rightarrow2$ a.s.

Let us study the tail distribution of $\tau_{1}.$ Since
\[
\mathbb{P}(\tau_{1}|X_{\tau_{1}}|>y|\xi_{\tau_{1}}=x)=\mathbb{P}(\tau
_{1}>y|\xi_{\tau_{1}}=x)=(1-|x|)^{y},
\]
by integration we obtain%
\[
\mathbb{P}(\tau_{1}>y)=\int_{-1}^{1}(1-|x|)^{y}|x|dx=2\int_{0}^{1}%
(1-x)^{y}xdx\sim2y^{-2}\text{ as }y\rightarrow\infty.
\]
Moreover, $\mathbb{E}(\tau_{k}X_{\tau_{k}})=0$ by symmetry. Also
\[
H(y)=\mathbb{E}(\tau_{1}^{2}I(\tau_{1}\leq y))\sim4\ln y\text{.}%
\]
Define a normalization satisfying $b_{n}^{2}\sim nH(b_{n}).$ In our case,
$b_{n}^{2}\sim4n\ln b_{n},$ implying that $b_{n}^{2}\sim2n\ln n.$

For each $n$, let $m_{n}$ be such that $T_{m_{n}}\leq n<T_{m_{n}+1}.$

We have the following representation
\begin{equation}
\sum_{k=1}^{n}X_{k}-\sum_{k=1}^{[n/2]}Y_{k}=(T_{0}-1)X_{0}+(\sum_{k=1}^{m_{n}%
}\tau_{k}X_{\tau_{k}}-\sum_{k=1}^{[n/2]}\tau_{k}X_{\tau_{k}})+\sum
_{k=T_{m_{n}+1}}^{n}X_{k}\text{,} \label{Rep}%
\end{equation}
where $Y_{k}=\tau_{k}X_{\tau_{k}}$ is a centered i.i.d. sequence in the domain
of attraction of a normal law. By the limit theorem for i.i.d. variables in
the domain of attraction of a stable law [see Feller (1971)] we obtain,%
\begin{equation}
\frac{\sum_{k=1}^{[n/2]}Y_{k}}{b_{[n/2]}}\Rightarrow N(0,1)\text{.}
\label{Fell}%
\end{equation}
By Theorem 4.1 from Billingsley (1968) the CLT for $(\sum_{k=1}^{n}%
X_{k})/b_{[n/2]}$ will follow from (\ref{Rep}) and (\ref{Fell}) provided we
show that the normalized quantity in the right-hand side of (\ref{Rep})
converges in probability to $0$. Clearly, because $b_{[n/2]}\rightarrow\infty$
we have
\[
\frac{(T_{0}-1)X_{0}}{b_{[n/2]}}\Rightarrow0.
\]
Also,%
\[
\mathbb{E}\frac{|\sum_{k=T_{m_{n}+1}}^{n}X_{k}|}{b_{[n/2]}}\leq\frac
{\mathbb{E}|\tau_{m_{n}+1}|}{b_{[n/2]}}=\frac{2}{b_{[n/2]}}\rightarrow0.
\]
Therefore it remains to study the middle term. Let $\delta>0.$%
\begin{gather*}
\mathbb{P}(|\sum_{k=1}^{m_{n}}Y_{k}-\sum_{k=1}^{[n/2]}Y_{k}|>\varepsilon
b_{[n/2]})\leq\mathbb{P}(|\frac{m_{n}}{n}-\frac{1}{2}|\geq\delta)\\
+\mathbb{P}(\max_{n/2-\delta n<l<n/2+\delta n}|\sum_{k=1}^{l}Y_{k}-\sum
_{k=1}^{[n/2]}Y_{k}|>\varepsilon b_{[n/2]})=I+II
\end{gather*}
Then, by the definition of $m_{n}$ and the law of large numbers for the i.i.d.
sequence $(\tau_{i})_{i\geq1}$ we know that: $m_{n}/n\rightarrow
1/\mathbb{E}(\tau_{1})=1/2$ a.s. Therefore the first term converges to $0$ for
every $\delta$ fixed as $n\rightarrow\infty$. As for the second term, by
stationarity and the fact that $Y_{k}$ are i.i.d.%
\[
II\leq2\text{ }\mathbb{P}(\max_{1\leq l\leq\lbrack\delta n]+1}|\sum_{k=1}%
^{l}Y_{k}|>\varepsilon b_{[n/2]}/2)
\]
and Theorem 1.1.5 in De la Pe\~{n}a and Gin\'{e} (1999),
\[
II\leq2\text{ }\mathbb{P}(\max_{1\leq l\leq\lbrack\delta n]+1}|\sum_{k=1}%
^{l}Y_{k}|>\varepsilon b_{[n/2]}/2)\leq18\text{ }\mathbb{P}(|\sum
_{k=1}^{[\delta n]+1}Y_{k}|>\varepsilon b_{[n/2]}/60)\text{.}%
\]
Then, by the central limit theorem in (\ref{Fell}) and the fact that
$b_{n}^{2}\sim2n\ln n,$ we have
\begin{align*}
\lim\sup_{n\rightarrow\infty}\mathbb{P}(|\sum_{k=1}^{[\delta n]+1}Y_{k}|  &
>\varepsilon b_{[n/2]}/60)=\lim\sup_{n\rightarrow\infty}\mathbb{P}(|\sum
_{k=1}^{[\delta n]+1}Y_{k}|/b_{[\delta n]}>\varepsilon b_{[n/2]}/60b_{[\delta
n]})\\
&  \leq\mathbb{P}(N(0,1)>\varepsilon\delta^{-1/2}/120)
\end{align*}
which converges to $0$ as $\delta\rightarrow0.$

It follows that%
\[
\frac{S_{n}}{b_{[n/2]}}\Rightarrow N(0,1)\text{.}%
\]
We recall that $\sigma_{n}^{2}=2n\log n=b_{n}^{2},$ implying that
\[
\frac{S_{n}}{\sigma_{n}}\Rightarrow N(0,\frac{1}{2})\text{.}%
\]
Consequently, because the chain is irreducible and aperiodic, by Corollary
\ref{corIA}, $W_{n}(t)\Rightarrow2^{-1/2}W(t)$. $\Diamond$

\bigskip

For a different example having this type of asymptotic behavior we cite Zhao
et al. (2010). Our Corollary \ref{corIA} will also provide a functional
central limit theorem for their example.

\bigskip

\textbf{Acknowledgement.} The authors would like to thank the referee for
carefully reading the manuscript and for numerous suggestions that improved
the presentation of this paper. The last author would like to thank Mikhail
Gordin for suggesting the use of  forward-backward martingale decomposition.

\end{document}